\documentclass{article}
\usepackage{latexsym, amsfonts}
\newcommand{\be}{\begin{equation}}
\newcommand{\ee}{\end{equation}}
\newcommand{\bea}{\begin{eqnarray}}
\newcommand{\eea}{\end{eqnarray}}
\newcommand{\bean}{\begin{eqnarray*}}
\newcommand{\eean}{\end{eqnarray*}}
\newcommand{\brray}{\begin{array}}
\newcommand{\erray}{\end{array}}
\newcommand{\ben}{\begin{equation}{nonumber}}
\newcommand{\een}{\end{equation}{nonumber}}

\newtheorem{dfn}{Definition}[section]
\newtheorem{thm}[dfn]{Theorem}
\newtheorem{lmma}[dfn]{Lemma}
\newtheorem{ppsn}[dfn]{Proposition}
\newtheorem{crlre}[dfn]{Corollary}
\newtheorem{xmpl}[dfn]{Example}
\newtheorem{rmrk}[dfn]{Remark}
\newcommand{\bdfn}{\begin{dfn}}
\newcommand{\bthm}{\begin{thm}}
\newcommand{\blmma}{\begin{lmma}}
\newcommand{\bppsn}{\begin{ppsn}}
\newcommand{\bcrlre}{\begin{crlre}}
\newcommand{\bxmpl}{\begin{xmpl}}
\newcommand{\brmrk}{\begin{rmrk}}
\newcommand{\edfn}{\end{dfn}}
\newcommand{\ethm}{\end{thm}}
\newcommand{\elmma}{\end{lmma}}
\newcommand{\eppsn}{\end{ppsn}}
\newcommand{\ecrlre}{\end{crlre}}
\newcommand{\exmpl}{\end{xmpl}}
\newcommand{\ermrk}{\end{rmrk}}

\newcommand{\IC}{\mathbb{C}}

\newcommand{\IR}{\mathbb{R}}

\newcommand{\cla}{{\cal A}}
\newcommand{\clb}{{\cal B}}
\newcommand{\clc}{{\cal C}}

\newcommand{\clf}{{\cal F}}
\newcommand{\clg}{{\cal G}}
\newcommand{\clh}{{\cal H}}

\newcommand{\clk}{{\cal K}}
\newcommand{\cll}{{\cal L}}

\newcommand{\clo}{{\cal O}}

\newcommand{\clq}{{\cal Q}}

\newcommand{\cls}{{\cal S}}
\newcommand{\clt}{{\cal T}}

\def\a*{{\cal A}_{h,*}}
\def\B{{\cal B}(h)}
\def\B1{{\cal B}_1(h)}
\def\b{{\cal B}^{\rm s.a.}(h)}
\def\b1{{\cal B}^{\rm s.a.}_1(h)}

\newcommand{\ot}{\otimes}

\newcommand{\raro}{\rightarrow}

\newcommand{\lgl}{\langle}
\newcommand{\rgl}{\rangle}

\def \qed {$\Box$}

\begin{document} 
\begin{center}
{\Large{\bf Quadratic independence of coordinate functions of certain homogeneous spaces and  action of compact quantum groups}}\\
{\large {\bf Debashish Goswami \footnote {The author gratefully acknowledges support from IUSSTF for Indo-US fellowship,  Dept. of Science and Technology  of Govt. of India for the Swarnajayanti fellowship and project and also the Indian National Science Academy.}}}\\
Indian Statistical Institute\\
203, B. T. Road, Kolkata 700108\\
Email: goswamid@isical.ac.in\\
\end{center}
\begin{center}
{\it Dedicated to the memory of Prof. Somesh Chandra Bagchi}
\end{center}
\begin{abstract}
Let  $G$ be one of the classical 
 compact, simple, centre-less, 
connected Lie groups or rank $n$ with a maximal torus $T$, the  Lie algebra $\clg$ and let $\{ E_i, F_i, H_i, i=1, \ldots, n \}$ be the standard set 
 of generators corresponding to a basis of the root system.
Consider 
 the adjoint-orbit space $M=\{ {\rm Ad}_g(H_1),~g \in G \}$, identified with the homogeneous space $G/L$ where $L=\{ g \in G:~{\rm Ad}_g(H_1)=H_1\}$. We prove that 
  the `coordinate functions' $f_i(g):=\lambda_i({\rm Ad}_g(H_1))$, $i=1, \ldots, n$, where $\{ \lambda_1, \ldots, \lambda_n\}$ is basis of $\clg^\prime$ are `quadratically 
   independent' in the sense that they do not satisfy 
   any nontrivial homogeneous quadratic relations among them. 
 Using this, it is proved that there is no genuine   compact quantum group which can act faithfully on   $C(M)$ such 
  that  the action  leaves invariant the linear span of the above coordinate functions. 
 
 As a corollary, it is also shown that any compact quantum  group having a faithful action on the noncommutative manifold obtained by  Rieffel deformation of $M$ 
 satisfying a similar `linearity' condition  must be  a Rieffel-Wang type deformation of some compact group.
\end{abstract}
\section{Introduction}
It is indeed a very important and interesting problem in the theory quantum groups and noncommutative geometry to study  `quantum symmetries' of various classical and quantum structures. Initiated by S. Wang who defined and studied quantum permutation groups of finite sets and quantum automorphism groups of finite dimensional algebras, such questions were taken up by a number of mathematicians including Banica, Bichon, and more recently by Goswami, Bhowmick and Skalski (see \cite{ban1}, \cite{ban2}, \cite{wang}, \cite{goswami}, \cite{qiso_comp}, \cite{qdisc} and also references therein).  Although there are many genuine (i.e. not of the form $C(G)$ for some group $G$) compact quantum groups which can faithfully act on the space of functions on finite sets, and these also give such actions on $C(X)$ for disconnected spaces with  at least $4$ components, no example of faithful action of a genuine compact quantum group on $C(X)$ for a connected compact space $X$ was known until recently, when H. Huang (\cite {
huchi}) gave a systematic way to construct  examples of such action of the quantum permutation groups. 
The work of Huang thus disproved a conjecture made by the author of the present paper about non-existence of faithful actions of genuine 
compact quantum groups on classical connected spaces. Let us also mention that using an algebraic example given by Etingof and Walton in 
 \cite{etingof_walton} it is possible to produce an example of faithful action by a finite dimensional compact quantum group 
  on the wedge sum of two copies of $([-1,1],0)$.

However, it is interesting to observe that none of the above examples  are smooth manifolds. 
On the other hand, it follows from the work of  Banica et al (\cite{banica_jyotish})  that most of known compact quantum groups, including the quantum permutation groups of Wang, can never act faithfully and isometrically on a connected compact Riemannian manifold. All these motivate us to modify the original conjecture to the following:\\
{\bf Conjecture:} It is not possible to have smooth faithful actions of genuine compact quantum groups on $C(M)$ when $M$ is a compact connected smooth manifold. 

 In this article, we provide some supporting evidence to this conjecture, by proving non-existence of nice (in some suitable sense to be described later) action of any genuine compact quantum group on 
 a large class of classical connected  manifolds which are homogeneous spaces of simple compact connected Lie groups. Using this, 
 we also prove that any compact quantum group acting faithfully and `linearly' in a certain sense on Rieffel-type deformation of $C(M)$  must be a Rieffel-Wang type deformation of 
 $C(K)$ for some classical group $K$ acting on $M$.
 
 We must mention that in a recent joint article with Joardar  (\cite{rigidity_general}), we have been able to prove this conjecture for isometric actions. 
 Nevertheless, the constructions and arguments of the present paper are quite  interesting in their own right, 
 for example, the quadratic independence of the  algebra of natural coordinate functions on a large class of homogeneous spaces of 
  a simple compact connected Lie group which is proved by a direct computation using the Serre relations only. 
 
Let us conclude this section with a brief remark on the physical relevance of the above fact. It has two implications: firstly, 
it implies that for a classical system with phase-space modeled on a compact connected manifold, the generalized notion of symmetries in terms of quantum groups coincides with 
the conventional notion, i.e. symmetries coming from group actions. This gives some kind of consistency of the philosophy of thinking quantum group actions as symmetries.
Secondly, it also allows us to describe all the  (quantum) symmetries of a physical model obtained by suitable deformation (at least for the Rieffel-type deformations) 
of a classical model with connected compact phase space, showing that such quantum symmetries are indeed deformations of the classical (group) symmetries of the original classical model.
 
\section{Preliminaries on quantum groups and their  actions}
\subsection{Basics of compact quantum groups}
A compact quantum group (CQG for short) is a  unital $C^*$ algebra $\cls$ with a coassociative coproduct 
(see \cite{woro1}) $\Delta$ from $\cls$ to $\cls \ot \cls$ (injective tensor product) 
such that each of the linear spans of $\Delta(\cls)(\cls \ot 1)$ and that of $\Delta(\cls) (1 \ot \cls)$ is norm-dense in $\cls \ot \cls$. 
From this condition, one can obtain a canonical dense unital $\ast$-subalgebra $\cls_0$ of $\cls$ on which linear maps $\kappa$ and 
$\epsilon$ (called the antipode and the counit respectively) making the above subalgebra a Hopf $\ast$ algebra. In fact, we shall always choose this dense 
 Hopf $\ast$-algebra to be the algebra generated by the `matrix coefficients' of the (finite dimensional) irreducible unitary representations (to be defined 
 shortly) of the CQG. The antipode is an anti-homomorphism and also satisfies $\kappa(a^*)=(\kappa^{-1}(a))^*$ for $a \in \cls_0$.
 
 It is known that there is a unique state $h$ on a CQG $\cls$ (called the Haar state) which is bi-invariant in the sense that $({\rm id} \ot h)\circ \Delta(a)=(h \ot {\rm id}) \circ \Delta(a)=h(a)1$ for all $a$. The Haar state need not be faithful in general, though it is always faithful on $\cls_0$ at least.

We say that a CQG $\cls$ (with a coproduct $\Delta$) (co)acts on a unital $C^*$ algebra $\clc$ if there is a unital $C^*$-homomorphism 
$\beta : \clc \raro \clc \ot \cls$ such that ${\rm Span}\{ \beta(\clc)(1 \ot \cls)\}$ is norm-dense in $\clc \ot \cls$, and it satisfies the coassociativity
 condition, i.e. $(\beta \ot {\rm id}) \circ \beta=({\rm id} \ot \Delta) \circ \beta$. It has been shown in \cite{podles} that 
 there is a unital dense $\ast$-subalgebra $\clc_0$ of $\clc$ such that $\beta$ maps $\clc_0$ into $\clc_0 \ot_{\rm alg} \cls_0$ (where $\cls_0$ is the dense 
 Hopf $\ast$-algebra mentioned before) and we also have 
   $({\rm id} \ot \epsilon)\circ
 \beta={\rm id}$ on $\clc_0$. In fact, this subalgebra $\clc_0$ comes from the canonical decomposition of $\clc$ into subspaces on each of 
which the action 
 $\beta$ is equivalent to an irreducible representation. More precisely, $\clc_0$ is the algebraic direct sum of finite dimensional vector spaces $\clc^\pi_i$, say, where $i$ runs over some index 
 set $J_i$, and $\pi$ runs over some subset (say $T$) of the set of (inequivalent) irreducible unitary representations of $\cls$, and the restriction of $\beta$ to $\clc^\pi_i$ is equivalent to the 
 representation $\pi$. Let $\{ a^{(\pi,i)}_j,j=1,...,d_\pi \}$ (where $d_\pi$ is the dimension of the representation $\pi$) be a basis of $\clc^\pi_i$ such that 
 $\beta(a^{(\pi,i)}_j)=\sum_k a^{(\pi,i)}_k \ot t^\pi_{jk}$, 
  for elements $t^\pi_{jk}$ of $\cls_0$. The elements $\{ t^\pi_{jk}, \pi \in T;~j,k=1,...,d_\pi \}$ are called the `matrix coefficients' of the action $\beta$.
 
 We say that the action $\beta$ is faithful if the $\ast$-subalgebra of $\cls$ generated by elements of the form $(\omega \ot {\rm id})(\beta(a))$, where $a \in \clc$, $\omega$ being a bounded linear functional on $\clc$, is norm-dense in $\cls$, or, equivalently, the $\ast$-algebra generated by the matrix coefficients is norm-dense in $\cls$.


A unitary representation of  a CQG $(\cls, \Delta)$ in a Hilbert space $\clh$ is given by a complex linear map  $U$ from the Hilbert space $\clh$ to the Hilbert $\cls$-module  $\clh \ot \cls$, which is isometric in the sense that $\lgl U\xi, U \eta \rgl_\cls=\lgl \xi, \eta \rgl$ for all $\xi, \eta \in \clh$ (where $\lgl \cdot, \cdot \rgl_\cls$ denotes the $\cls$-valued inner product) and the $\cls$-linear span of the range of $U$ is dense in $\clh \in \cls$.  
 There is an equivalent description of the unitary representation given by the $\cls$-linear unitary $\tilde{U} \in \cll(\clh \ot \cls)$ defined by $\tilde{U}(\xi \ot b)=U(\xi) b$, 
for $\xi \in \clh, b \in \cls,$ satisfying $({\rm id} \ot \Delta) (\tilde{U})=\tilde{U}^{12}\tilde{U}^{13}$. Here, $\tilde{U}^{ij}$ is the standard leg-numbering notation used in the theory of quantum groups, i.e. $\tilde{U}^{12}=  \tilde{U} \ot {\rm Id}_{\cls}$ and $\tilde{U}^{13}:=\sigma_{23} \circ \tilde{U}^{12} \circ \sigma_{23},$ $\sigma_{23}$ being the map which flips the second and third tensor components of $\clh \ot \cls \ot \cls$.
We  also recall that a unitary representation $U$ is irreducible if there is no nontrivial closed subspace $\clk$ of $\clh$ which is invariant under $U$, i.e. $U\clk \subseteq \clk \ot \cls$.

\section{Quadratic independence and nonexistence of genuine quantum group action}
Let $V$ be a finite dimensional subspace of a (possibly infinite dimensional) commutative algebra say $\cla$ over any field $F$. Let $V^{(2)}$ denote the linear span of elements of the form 
$vv^\prime,$ with $v,v^\prime \in V$.
 Let $V\ot_{\rm sym}V$ denote the symmetric tensor product of $V$ with itself, which is the subspace $V \ot V$ spanned by elements of the form $v \ot v$, $v \in V$. Equivalently, it is the 
   vector space obtained by taking quotient $V \ot V$ by the subspace $\clf$ spanned by vectors of the form $v \ot w-w \ot v$. As the algebra $\cla$ is commutative, the linear 
    map $v \ot w \mapsto vw=wv$ from $V \ot V$
    to $V^{(2)} \subseteq \cla$ annihilates the subspace $\clf$, hence induces a linear map from the quotient, i.e. $V \ot_{\rm sym} V$ onto $V^{(2)}$. We call $V$ to be {\it quadratically independent}
     if this map is one-to-one, i.e. a linear isomorphism. 
   We also call the dimension of $V^{(2)}$ the quadratic dimension of $V$.  
\blmma
The following are equivalent:\\
(i) $V$ is quadratically independent.\\
(ii) The quadratic dimension of $V$ equals $\frac{n(n+1)}{2}$, where $n={\rm dim}(V)$.\\
(iii) For every basis $\{x_1,...,x_n\}$ of $V$, the set $\{ x_ix_j,~1\leq i\leq j \leq n\}$ is linearly independent.
\elmma
{\it Proof:}\\
To see the equivalence of (i) and (ii), it is enough to note that there is a surjective linear map from the finite-dimensional space $V \ot_{\rm sym} V$ 
 to $V^{(2)}$, hence these two spaces are isomorphic if and only if dimension of $V^{(2)}$ equals that of  $V \ot_{\rm sym} V$, which is nothing but 
  $\frac{n(n+1)}{2}$. It is also clear that (iii) implies (ii). To prove the other way, observe that for any basis $\{ x_1,...,x_n \}$ of $V$, the set $E=\{ x_ix_j, i\leq j \leq n\}$ spans $V^{(2)}$, hence a subset of $E$ will give a basis of $V^{(2)}$. But $E$ has cardinality not greater than $n(n+1)/2$, which is by (i) the dimension of $V^{(2)}$. Therefore, $E$ must be a basis itself.
\qed

We remark that if the underlying field is $\IR$, then quadratic independence of $V$ is clearly equivalent to the quadratic independence of the complexification $V_{\IC}$, as $(V^{(2)})_{\IC}\cong (V_{\IC})^{(2)}$.

We end this section with the following crucial implication of quadratic independence in the context of quantum group actions. 
 We refer the reader to \cite{rigidity_general} (Lemma 10.1, see also Theorem 2.2 of \cite{qiso_comp}) for a proof.
\bthm
\label{no_quantum}
Let $\cla$ be a unital commutative $C^*$ algebra and $x_1,...,x_n$ be self-adjoint elements of $\cla$ such that ${\rm Span}\{ x_1,...,x_n\}$ is quadratically independent and $\cla$ is the unital $C^*$ algebra generated by $\{x_1,...,x_n\}$. Let $\clq$ be a compact quantum group acting faithfully on $\cla$, such that the action maps  
${\rm Span}\{ x_1,...,x_n\}$ into itself. Then $\clq$ must be commutative as a $C^*$ algebra, i.e. $\clq \cong C(K)$ for some compact group $K$.
\ethm
\section{Some facts about simple Lie groups and Lie algebras}
We collect a few standard facts about simple Lie groups and Lie algebras. Most of these materials are taken from \cite{helga}, \cite{hump}. 
Let $G$ be a compact, simple, connected Lie group with trivial center, and let $\clg$ be its (real) Lie algebra.  Consider the complexification $\clg_{\IC}$ of $\clg$, and the corresponding simple 
 Lie group $G_{\IC}$, which has the (complex) Lie algebra $\clg_{\IC}$ (see \cite{helga} pages 178-182).  For $X \in \clg$, let ${\rm ad}_X(\cdot)=[X. \cdot] : \clg \raro \clg$ and for $g \in G$ 
  let ${\rm Ad}_g$ denote the adjoint action of $G$ at $g$, which is obtained as the differential of the map $G \ni h \mapsto ghg^{-1} \in G$ at $h=e$ (identity element of $G$).
   The adjoint representation of $G$ (see page 127 of \cite{helga}), i.e. $g \mapsto {\rm Ad}_g$ identifies $G$ with a matrix group acting on $\clg$, i.e. as subgroup of $GL(\clg)$. 
   We refer the reader to page 104, \cite{helga} for the definition and properties of the exponential map ${\rm exp}$ and 
 shall use the abbreviation $\beta_t(X)={\rm Ad}_{{\rm exp}(tX)}={\rm exp}(t {\rm ad}_X)$ for $X \in \clg_\IC$. We also note that the adjoint representation is an irreducible representation 
  of $G$ as $G$ has a trivial center.

  Suppose that $\clg$ has rank 
 $n$. Let $T$ be a maximal torus with the corresponding maximal abelian Lie subalgebra (Cartan algebra) $\clt$.  
 From the general structure theory of simple Lie algebras (see, e.g page 96 of \cite{hump}), $\clg_{\IC}$ is  generated (as Lie algebra)
 by elements  $\{ E_i,F_i,H_i:~i=1,2,\ldots,n\}$, satisfying the following relations (`Serre relations'):\\
$$ [E_i,F_i]=H_i,~~[E_j,H_i]=-a_{ij}E_j,~~[F_j,H_i]=a_{ij}F_j,$$
$$[E_i,F_j]=0~\forall i\neq j,~~[H_i,H_j]=0~\forall i,j,$$
$$ {\rm ad}_{E_i}^{1-a_{ij}}(E_j)=0,~~{\rm ad}_{F_i}^{1-a_{ij}}(F_j)=0,~~i \neq j.$$ Here, $\{ H_1,\ldots, H_n\}$ is a set of generators of $\clt$ and  $(( a_{ij} ))$ is the so-called 
the  Cartan matrix (page 55, \cite{hump}), with $a_{ii}=2,$ $a_{ij} \leq 0$ for $i \neq j$.
Thus, a spanning set (as vector space) of the Lie algebra $\clg_\IC$ consists of  $E_i, F_i, H_i$, $i=1,...,n$, and those elements 
of the form ${\rm ad}_{E_{i_1}}^{m_1}...{\rm ad}_{E_{i_k}}^{m_k}(E_{i_{k+1}}),$ ${\rm ad}_{F_{i_1}}^{n_1}...{\rm ad}_{F_{i_l}}^{n_l}(F_{i_{l+1}})$, $1 \leq i \leq n,k,l \leq n-1,i_1<...<i_{k+1}$ 
which are nonzero. Let us denote by $\hat{\clb}$ this spanning set. 

Given an element $X$ of $\clg$, consider the adjoint-orbit space  $\clo_X=\{ {\rm Ad}_g(X):~g \in G \}$  which is a closed subset of $\clg$. It is easy to see that $\clo_X$ 
can be identified 
 as a manifold with the homogeneous space $G/L$ where $L=\{ g:~{\rm Ad}_g(X)=X\}$. Indeed, it is clear from the definition that $L$ is a closed subgroup, hence a Lie subgroup. 
 Thus, its Lie algebra, say $\clg_X$, will
  be a Lie subalgebra of $\clg$. In fact, $Y \in \clg_X$ if and only if ${\rm exp}(tY) \in L$  for all $t \in \IR$, i.e. $\beta_t(Y)(X)=0$ for all $t$. This is equivalent to $[Y,X]=0$. 
   In other words, $\clg_X=\{ Y \in \clg:~ [Y,X]=0\}$. 
  
  In this paper, we shall consider $G$ to be one of the classical groups (i.e. special unitary, special orthogonal or symplectic groups) and $X=H_1$ 
  (equivalently, $H_i$ for some fixed $i$).\\
  
    The following fact is easily verifiable from the explicit description of the group of inner automorphisms, in particular the Weyl groups for the classical simple Lie groups:\\
   {\it Fact:} For each $i$, there is $g_i \in G_\IC$ such that ${\rm Ad}_{g_i}(H_1)=H_i$.
   
 \brmrk
 We do not know whether the above fact is true for exceptional Lie algebras as well. If so, the above fact and hence all the results of this paper will hold for an 
 arbitrary simple compact 
  connected centreless Lie group.
 \ermrk

 Now, as $\clg$ is a simple Lie algebra,  its Cartan matrix is indecomposable (page 73, \cite{hump}). That is, for every $i$ and $j$, we can find a `path',
 say $i_1=i,i_2,...,i_k=j$, such that $a_{i_l i_{l+1}} \neq 0$ for each $l=1,...,k-1$. Unless otherwise mentioned, we shall always assume that no two induces
 in a `path' are chosen equal, and let us 
 say that the length of such a path is $k-1$. Denote by $\Gamma^k_m$ the subset of $i$ such that we can find a `path' of length at most $m$ from $k$ to $i$ and 
 let $\hat{\clb}^k_m$ be  the subset of 
 $\hat{\clb}$ consisting of those $E_i,F_i, H_i$ with $i \in \Gamma^k_m$  and also 
 $[E_{i_1},...,E_{i_k}],[F_{j_1},...,F_{j_l}]$ with all $i_p$ and   $j_q$  belonging to $\Gamma^k_m$. Choose a linearly independent subset $\clb$ of $\hat{\clb}$ as follows. 
Start with all $E_i, F_i, H_i$'s in $\Gamma^1_1$ (which are clearly linearly independent), then extend this collection  to a linearly independent set, say $\clb_1$, 
which has the same span as that of ${\hat{\clb}}_1$ by adding appropriate elements from ${\hat{\clb}}^1_1$. Then add further suitable elements from 
${{\hat{\clb}}^1}_2 \bigcap ({{\hat{\clb}}^1}_1)^c$ 
to get a linearly independent set $\clb_2$ with the same span as that of ${\hat{\clb}}^1_2$. Go on like this to get $\clb=\bigcup_m \clb_m$, which is a basis of $\clg$, 
with $\clb_m=\clb \bigcap {\hat{\clb}}^1_m$. Here $\clb_m$ is a basis for ${\rm Span}({\hat{\clb}}^1_m)$ for all $m$.
Note that $\clb$ is actually a union of finitely many $\clb_m$'s due to the Serre relations. 
\section{Main results}  
\subsection{Nonexistence of linear quantum actions on the adjoint-orbit space }

The real $C^*$ algebra of real-valued continuous functions $C(M)_\IR$ of the adjoint-orbit space $M=G/L$ can be identified with real $C^*$-sub algebra of $C(G)_\IR$ 
generated by the functions of the 
form $f_\lambda(g)=\lambda({\rm Ad}_g(H_1)),$ where $\lambda$ varies over the dual ${\clg}^\prime$ of $\clg$.  Let $V$ be 
the real vector space consisting of functions of the form $f_\lambda$, 
with $\lambda \in \clg^\prime$. We can complexify all the vector spaces involved, and observe that $V_\IC$ will generate $C(M)$ 
(continuous complex valued functions) as a (complex) $C^*$ algebra.
The real-valued functions (i.e. self-adjoint as elements of $C(M) \subset C(G)$) $\{ f_i \equiv f_{\lambda_i},~i=1, \ldots,n \}$, where 
 $\{ \lambda_1, \ldots, \lambda_n\}$ is a basis of ${\clg}^\prime$, can be interpreted as `coordinates' of the compact homogeneous manifold $M$. 

Note that  the canonical action of $G$ on the homogeneous space $M=G/L$ leaves invariant the subspace $V$ mentioned above, i.e.
it is `linear' in the sense discussed before. It is thus natural to consider linear actions by quantum groups as well. In fact, very often the natural invariant 
Laplacian on $M$ given by the action of the Casimir will have $V$ as an eigenstate, hence any isometric action by a compact quantum group must leave $V$ invariant and is determined by the restriction of the action to $V$. 

However, as we'll see below, there is no genuine quantum group action on $M$ which is linear.


We claim the following.
\bthm
\label{main}
The vector space $V$ is quadratically independent.
\ethm
 The proof will be divided into several steps. First of all, we need to pass  to the complex Lie group $G_\IC$.
 \blmma
 Consider the natural extension of adjoint action of $G_\IC$ on $\clg_\IC$,  and let $\tilde{V_\IC}$ be the linear subspace spanned by functions of the form $\tilde{f_\lambda}(g):=\lambda({\rm Ad}_g(H_1))$ ($g \in G_\IC$) in $C(G_\IC)$. Then the quadratic independence of $V$ is equivalent to quadratic independence of $\tilde{V_\IC}$.
 \elmma
 {\it Proof:}\\
 We already noted that quadratic independence of $V$ and $V_\IC$ are equivalent. Moreover, as $\tilde{f_\lambda}$'s are extensions of $f_\lambda$'s, 
 it is clear that quadratic independence of $V_\IC$ will imply that of $\tilde{V_\IC}$. To prove the converse, assume that $\tilde{V_\IC}$ is quadratically independent. 
 We have to show that $V_\IC$ is quadratically independent. Choose a basis $f_1,...,f_m$ of $V_\IC$; clearly $\tilde{f}_1,...,\tilde{f}_m$ are linearly independent and 
 hence forms a basis of $\tilde{V_\IC}$. If possible, suppose that $c_{ij}, i\leq j$ are complex numbers such that $\sum_{ij} c_{ij} f_if_j=0$. This means for every 
 element $X$ of $\clg$, and $t \in \IR$, we have $\sum_{ij} c_{ij} f_i({\rm exp}(tX))f_j({\rm exp}(tX))=0$. Now, it is known that $z \mapsto {\rm exp}(zX)$ is holomorphic, 
 and from this, it can be seen that the above expression with $t$ replaced by $z \in \IC$ is holomorphic, so that it must be identically zero for all $z \in \IC$.
  But then, as $\clg_\IC$ is the Lie algebra 
of the connected Lie group $G_\IC$, we conclude that $\sum_{ij} c_{ij} \tilde{f_i}(g)\tilde{f_j}(g)=0$ for all $g \in G_\IC$. By assumption, 
this implies $c_{ij}=0$ for all $i,j$, i.e. $V_\IC$ is proved to be quadratically independent.
 \qed
 
In view of this lemma, we can prove the theorem by showing quadratic independence of $\tilde{V_\IC}$. Let us first show  that the  vector space of functions $\tilde{f_\lambda}$ has the same dimension as the dimension of $\clg_\IC$. 
\blmma
 The map $\clg^\prime \ni \lambda \mapsto \tilde{f}_\lambda \in \tilde{V}_\IC$ is an isomorphism of vector spaces.
 
 \elmma
 {\it Proof:}\\
 We need only to show injectivity of the map $\clg^\prime_\IC \ni \lambda \mapsto \tilde{f}_\lambda$. Let  $\lambda$ be such that 
  $\lambda({\rm Ad}_g(H_1)=0$ Let $\clh$ be the subset of $\clg_\IC$ consisting of all $Y$  such that $\lambda({\rm Ad}_g(Y))=0$ for all $g$. $\clh$ is nonzero as 
   it contains $H_1$ by assumption.
    Now, $\clh$ is clearly a subspace  and moreover, it is a Lie ideal in $\clg_\IC$. Indeed, if $Y \in \clh$ and $Z$ is arbitrary element of $\clg_\IC$, considering $h \in G_\IC$, 
    $t \in \IR$ and $g=h{\rm exp}(tZ)$
    we see $\lambda({\rm Ad}_{h{\rm exp}(tZ)}(Y))=0$. Differentiating at $t=0$ we get $\lambda({\rm Ad}_h([Z,Y]))=0$ for all $h \in G_\IC$, hence $[Z,Y] \in \clh$. But $G_\IC$ being simple 
     Lie algebra, any nonzero Lie ideal must be $\clg_\IC$ itself. This implies $\lambda({\rm Ad}_g(Y))=0$ for all $Y$ and for all $g$, so in particular $\lambda(Y)=0$ for all $Y$, i.e. $\lambda=0$.
     \qed

 Using this lemma, we can identify a functional $\lambda$ on $\clg_\IC$ with the function on $G$ given by $g \mapsto \lambda({\rm Ad}_g(H_1))$, 
 so $\lambda \ot \eta$ (for $\lambda, \eta \in 
  \clg^\prime$),  can be identified with $\tilde{f}_\lambda(g)f_\eta(g)=(\lambda \ot \eta)({\rm Ad}_g(H_1),{\rm Ad}_g(H_1))$. 
  More generally, $c \in \clg_\IC^\prime \ot \clg_\IC^\prime$ is  identified
   with $g \mapsto \hat{c}(g):=c({\rm Ad}_g(H_1),{\rm Ad}_g(H_1))$. Thus, the image of the symmetric tensor product of $V_\IC$ with itself in the algebra $C(G)$ 
    consists of the functions of the form $\hat{c}$ as above with $c \in \clg_\IC^\prime \ot_{\rm sym} \clg_\IC^\prime$. 
 By definition of quadratic 
  independence, it suffices to show that $c \mapsto \hat{c}$ is one-to-one.
 

{\it Proof of quadratic independence}:\\
Consider a bilinear symmetric functional $c$ on $\clg_\IC \ot \clg_\IC$, i.e. element of $\clg^\prime_\IC \ot_{\rm sym} \clg^\prime_\IC$ and let $\hat{c}$ 
 denote its canonical image in $C(M)$ given by $\hat{c}(g)=c({\rm Ad}_g(H_1),{\rm Ad}_g(H_1))$ as discussed above.
 We'll usually write $c(X,Y)$ as $c_{(X,Y)}$.  We want to prove the injectivity of  $c \mapsto \hat{c}$. To this end, suppose $\hat{c}=0$, i.e. $\hat{c}(g)=0$ for all $g$. 
 We have to show that each $c_{(X,Y)}=0$.

First, let us note the following relations:\\
\be
\label{e}
\beta_t(E_i)(H_1)=H_1-ta_{1i}E_i,
\ee
\be
\label{f}
\beta_t(F_i)(H_1)=H_1+ta_{1i}F_i,
\ee
\be
\label{[ee]}
\beta_t([E_i, [E_j,...,E_k]...])(H_1)=H_1-(a_{1i}+a_{1j}+...+a_{1k})t[E_i,[E_j,...,E_k]...],
\ee
\be
\label{[ff]}
\beta_t([F_i,[F_j,...,F_k]...])(H_1)=H_1+(a_{1i}+a_{1j}+...+a_{1k})t[F_i,[F_j,...,F_k]...],
\ee
\be 
\label{h}
\beta_s(E_i)\beta_t(F_i)(H_1)=H_1-sa_{1i}E_i+ta_{1i}F_i+st a_{1i} H_i- s^2t a_{1i}E_i.
\ee

 For $X=[E_{i_1},[ \ldots,[E_{i_l},\ldots]]$ (respectively $[F_{i_1},[ \ldots,[F_{i_l},\ldots]]$) we denote by $a_{kX}$ the sum $a_{ki_1}+\ldots+a_{k i_l} +\ldots$ ($-a_{ki_1}-\ldots$ respectively),
  so that we have \be \label{zzzz} \beta_t(X)(H_k)=H_k-a_{kX}X. \ee Also, note that $a_{k [X,Y]}=a_{kX}+a_{kY}$.
 
 The strategy is similar to the proof of the previous lemma.  We first give explicit computations in certain simpler yet illustrative special cases and then give the general 
  inductive scheme.
 Consider
 
\be \label{maineqn2}\hat{c}({\rm Ad}_g(H_1), {\rm Ad}_g(H_1))=0, g=\beta_{t_1}(X_1) \ldots \beta_{t_k}(X_k).\ee

{\it Step 1}: $c_{(H_i,X)}=0$ for all $i$ and for $X=H_i,E_i,F_i$.\\
Consider $g_i$ such that ${\rm Ad}_{g_i}(H_1)=H_i$ and take $g=\beta_t(E_i)g_i$ in (\ref{maineqn2}).
Equating coefficients of $1$ and $t$ to $0$ and noting $a_{ii}=2$ we have $c_{(H_i,H_i)}=0$ and $c_{(H_i,E_i)}=0$ respectively. Similarly, we get $c_{(H_i,F_i)}=0$ too.\vspace{2mm}\\
{\it Step 2}: $c_{(E_i,F_i)}=0 \forall i$ and for $k$ and $i$ be such that $i \in \Gamma^k_1$, we have  $c_{(H_k,H_i)}=0.$\\
Using the analogue of (\ref{h}) with $k$ replacing $1$ in (\ref{maineqn2}) with
  the choice $g=\beta_s(E_i)\beta_t(F_i) g_k$ (where ${\rm Ad}_{g_k}(H_1)=H_k$) and comparing coefficients of $s^2t^2$ one 
gets \be \label{1} -2 a_{ki}^2c_{(E_i,F_i)}+a_{ki}^2c_{(H_i,H_i)}=0,~i.e.~c_{(H_i,H_i)}=2c_{(E_i,F_i)}.\ee As
 $c_{(H_i,H_i)}=0$ from Step 1, we conclude $c_{(E_i, F_i)}=0$. 
Similarly, from coefficient of $st$ we get \be \label{2} c_{(H_k,H_i)}=a_{ki}c_{(E_i,F_i)}.\ee This gives $c_{(H_k,H_i)}=0$.\vspace{2mm}\\
{\it Step 3}:\\Let $X,Y $ be of the form $[E_{i_1},[ \ldots,[E_{i_m}],\ldots]]$ or $[F_{j_1},[\ldots,F_{j_n}]\ldots]$, i.e. Lie Brackets of $E$'s or of $F$'s only. 
Then we have the following:\\
(a) $c_{(H_k,X)}=0$ if $a_{kX}$ is nonzero.\\
(b) If $a_{kX},a_{kY}$ are nonzero, then $c_{(X,Y)}$ is a nonzero multiple of $c_{(H_k,[X,Y])}$ and hence it is $0$ if and only if $c_{(H_k,[X,Y])}=0$.\\

We prove (a) by considering coefficient of $t$ in (\ref{maineqn2}) using (\ref{zzzz}). For (b), consider $\beta_s(X)\beta_t(Y)(H_k)$ and 
look at the coefficient of $st$ in (\ref{maineqn2}).\vspace{2mm}\\
{\it Step 4}: The conclusions of Steps 1,2,3 hold even if we replace $X$, $Y$, $E_i$, $F_i$, $H_i$ etc. by ${\rm Ad}_g(\cdot)$ of them. For example, 
 we have $c_{({\rm Ad}_g(X), {\rm Ad}_g(H_k))}=0$ if $a_{kX}$ is nonzero, $c_{({\rm Ad}_g(E_i),{\rm Ad}_g(F_i))}=0$ etc. \\
 The proofs are essentially the same as the previous steps. For example, for proving the analogue of Step 3 (a)for a fixed $g_0$, we need to consider the coefficient of 
  $t$ in (\ref{maineqn2}) using $g=\beta_t({\rm Ad}_{g_0}(X))g_0g_k$.\vspace{2mm}\\
    {\it Step 5}: For any $X$ which is either a Lie Bracket of only $E_i$'s or of $F_i$'s we have\\
    (a) $c_{({\rm Ad}_g(X), {\rm Ad}_g(X))}=0$  $\forall g$;\\
    (b) $c_{([Y,{\rm Ad}_g(X)],{\rm Ad}_g(X))}=0$  $ \forall Y \in \clg$, $\forall g$.\\
        We prove (a) by noting that given any nonzero $X$ as above, there is some $i$ for which $a_{iX}$ is nonzero, because otherwise $X$ will commute with all $H_i$, hence must be an element of 
     the Cartan algebra. Now (a) follows from Step 4. The part (b) is obtained by replacing $g$ in (a) by $\beta_t(Y) g$ and then differentiating at 
     $t=0$ as well as using the fact that $c(\cdot, \cdot)$
      is symmetric in its arguments.\vspace{2mm}\\
{\it Step 6}: $c_{(H_k,Z)}=0$ for all $k$ and  $Z$.\\
This needs an inductive argument. Without loss of generality, suppose that $Z=[E_{i_1},[E_{i_2}, \ldots, ]]$, as the cases with $F_i$'s can be treated 
 similarly. First, if none of the $i_l$ are equal to $k$ and also at least one of them is in $\Gamma^k_1$ clearly $a_{kZ}$ is nonzero, so that $c_{(H_k,Z)}=0$. 
 Next, 
 let us  give 
  a simple special case where $a_{kZ}$ is nonzero by just adding one $E_k$ in the chain. 
    Consider $Z^\prime=[E_k,Z]$ and assume $a_{kZ}=-2$ so that $a_{kZ^\prime}=0$. However, we still want to show $c_{(H_k,Z^\prime)}=0$. We have by Step 5,
  taking $Y={\rm Ad}_g(E_k)$ and $X=Z$, 
 the following:
   $$c_{({\rm Ad}_g(Z^\prime),E_k)}=0 \forall g.$$ Choosing $g=\beta_t(F_k)$ and differentiating at $t=0$ gives $c_{(Z^\prime, H_k)}=c_{([F_k,Z^\prime],E_k)}.$ But as $Z$ has no $E_k$ in its 
    expression we have $[F_k,Z]=0$, hence $[F_k,Z^\prime]=-a_{kZ}Z=2Z$. Moreover, we already know $c_{(Z,H_k)}=0$, which proves $c_{(H_k,Z^\prime)}=0$ too.
    
    Now we come to the inductive scheme for proving Step 6.  Given any $Z$ which is a Lie Bracket of only $E_i$'s or only $F_j$'s, let us call it is at a distance $m$ 
    ($m \geq 1$) from $H_k$ if $m$ is the smallest positive integer for which 
     we can find $j_1=k, j_2, \ldots,j_{m}$  so that $a_{j_l j_{l+1}}$ and $a_{j_{m}Z}$ are  nonzero for all $l$. We make a similar definition for $Z=H_i$
      replacing the condition  $a_{j_{m}Z} \ne 0$ by $a_{j_{m}i} \ne 0$. We make the 
  induction hypothesis that $c_{(H_k,Z)}=0$ for any $Z$ which is at a distance less than or equal to $m$. For $m=1$ we already know it is true.  
  
  Assume the induction hypothesis for $m=N$ and let us first prove $c_{(H_k,H_j)}=0$ for 
 $k,j$ with $j \in \Gamma^k_{N+1}$.   Without loss of generality we can assume that $j$ does not belong to $\Gamma^k_N$ (
    there is nothing to prove otherwise), and let  $k=j_1,...,j_{N+2}=j$  being the corresponding minimal path of length $N+1$ from $k$ to $j$ (the $j_l$'s are distinct). 
Consider
     the following function $\xi\equiv \xi(w, w^\prime,s_{N+1}, t_{N+1}, \ldots, s_2, t_2)$ of $2N+2$ real variables, given by 
$$\xi=\beta_w(E_j)\beta_{w^\prime}(F_j) \beta_{s_{N+1}}(E_{j_{N+1}})\beta_{t_{N+1}}(F_{j_{N+1}})\ldots \beta_{s_2}(E_{j_2})\beta_{t_2}(F_{j_2})({\rm Ad}_{g_k}(H_1)),$$ 
where $g_k$ is such that ${\rm Ad}_{g_k}(H_1)=H_k$.
 We have $\hat{c}(\xi, \xi)=0$ and want to equate coefficients of different powers of the above real variables to $0$. For a $2N+2$-tuple
 $\underline{l}=(l_{N+2},l^\prime_{N+2}, l_{N+1},
  l^\prime_{N+1}, \ldots, l_2, l^\prime_2)$ of nonnegative integers, we denote by $\xi_{\underline{l}}$ 
  the coefficient of $w^{l_{N+2}}(w^\prime)^{l^\prime_{N+2}}s_{N+1}^{l_{N+1}} \ldots (s_2^\prime)^{l^\prime_2}$
   in the power series expansion of $\xi(\cdot)$. The coefficient of $ww^\prime s_{N+1}t_{N+1} \ldots s_2t_2$ in $\hat{c}(\xi, \xi)$ is 
   the sum of all $c(\xi_{\underline{l}}, \xi_{\underline{m}})$
     where $\underline{l}, \underline{m}$ are such that their coordinate-wise sum is $(1,1,\ldots, 1).$  
A moment's reflection would show that in this expression the coefficient of $s_2t_2...s_{N+1}t_{N+1}ww^\prime$, i.e. 
 $\xi_{(1,1,\ldots, 1)}$  is a nonzero  multiple of $H_j$.  Note that by minimality of the path chosen, $a_{j_pj_{p+r}}=0$ for any $p$ and any $r>1$, and so $E_{j_p}$ and $F_{j_q}$ commute
 whenever $|p-q|\geq 2$. 
 Moreover, observe the following:\\
 (i) If $l_2=l_2^\prime=0$ but $\underline{l}$ is not $(0,0, \ldots, 0)$, then $\xi_{\underline{l}}=0$.\\
 (ii) In fact, except for the case $(0,0,\ldots, 0)$, if any two consecutive entries of the finite sequence $\underline{l}$ are zero 
  then $\xi_{\underline{l}}=0$ too.\\
  Indeed, as $a_{j_p j_q}=0$ for $q >p+1$, we have $\beta_{t_{j_q}}(F_{j_q})(H_k)=\beta_{s_{j_q}}(E_{j_q})(H_k)=H_k$ for all $q \geq 3$,  hence 
   $\xi(w,w^\prime,s_{N+1}, t_{N+1}, \ldots, s_3, t_3,0,0)$ is identically equal to $H_k$, which proves (i).  Proof of (ii) is similar. 
   Now, (i),(ii) imply that the only 
   choices of $(\underline{l},\underline{m})$ which give 
   possibly nonzero contribution are $((0,0,\ldots,0), (1,1,\ldots,1))$ and $((1,0,1,0,\ldots,)(0,1,0,1,\ldots))$, so that 
   the coefficient of $s_2t_2...s_{N+1}t_{N+1}ww^\prime$ is a linear combination of $c_{(H_k,H_j)}$ and $c_{(X,Y)}$ 
   where $X=[E_{j_{N+1}},[E_{j_N},\ldots, E_{j_2}]\ldots]$, $Y=[F_{j_{N+1}},[F_{j_N},\ldots, F_{j_2}]\ldots]$. 
    But $a_{j_2X}$ and $a_{j_2Y}$ are nonzero, because $a_{j_2j_3}$ is nonzero and $a_{j_2j_p}=0$ for all $p>3$. 
     Moreover, $[X,Y]$ is a linear combination of $H_{j_p}$'s, $p=2,3,\ldots$ so that by the induction hypothesis, 
      $c_{(H_{j_2},[X,Y])}=0$. This implies $c_{(X,Y)}=0$ by Step 3(b). This proves $c_{(H_k,H_j)}=0$, as we have already noted that $c_{(H_k,H_j)}$ occurs with 
       a nonzero multiple in the expression.
     
     A similar argument considering $\beta_w(X) \beta_{s_{N+1}}(E_{j_{N+1}})\beta_{t_{N+1}}(F_{j_{N+1}})\ldots \beta_{s_2}(E_{j_2})\beta_{t_2}(F_{j_2})({\rm Ad}_{g_k}(H_1))$ will prove 
      $c_{(H_k, X)}=0$ for any $X$ which is Lie brackets of only $E$'s or of $F$'s, and at a distance $N+1$ from $k$, completing the proof of the inductive step, 
       hence also that of Step 6.\vspace{2mm}\\

  {\it Step 7}: The general case\vspace{2mm}\\
  Consider any two elements $X,Y$ of the Lie algebra. If $X=H_k$ for some $k$, we have already proved $c_{(X,Y)}=0$. So, assume $X$ is either Lie bracket of $E$'s or that of $F$'s. 
  Now, note that by a slight modification of Step 6, we get $c_{({\rm Ad}_g(H_k),{\rm Ad}_g(Y))}=0$ for all $g$, $k$ and $Y$. As before, this will give us $c_{([X,H_k],Y)}=-c_{(H_k,[X,Y])}
  =0$, i.e. $a_{kX}c_{(X,Y)}=0$. Next, choose $k$ such that $a_{kX}$ is nonzero, which proves $c_{(X,Y)}=0$.

\qed\vspace{2mm}\\

From the above theorem and Theorem \ref{no_quantum}, we immediately conclude the following:
\bthm
\label{main_no_action}
Let $G$ be a classical compact, connected, centreless  simple Lie group of rank $n$  as above, with $\{ H_1, \ldots, H_n\}$ a basis of the Cartan subalgebra 
 corresponding to simple roots. Let $L=\{ g \in G:~{\rm Ad}_g(H_1)=H_1 \}$ and $M=G/L$ be the corresponding homogeneous space. One can identify $C(M)$ with the subalgebra of 
  $C(G)$ consisting of $f$ satisfying $f(gh)=f(g)$ for all $h \in L$ and for all $g \in G$. Let 
   $V \subset C(M)$ be the subspace of functions of the form $g \mapsto \lambda({\rm Ad}_g(H_1)), \lambda \in \clg^\prime$.  Then any compact quantum group acting faithfully on $C(M)$  
such that $V$ is left invariant by the action must be commutative as a $C^*$ algebra.
\ethm

\brmrk
If one tries to generalize the techniques of this paper to semi-simple (but not necessarily simple) Lie algebras, one faces the following problem. Suppose that the Lie algebra 
 is a direct sum of $k$ copies of simple classical Lie algebras. Then the natural analogue of the space $V$ considered in this paper is the tensor product of the corresponding 
  $V_i$'s, $i=1,\ldots,k$ coming from the constituent simple Lie algebras. But this tensor product is not going to be quadratically independent for $k >1$. However, the conclusion 
   of Theorem \ref{main_no_action} may  still follow from the more general result of \cite{rigidity_general} as  already remarked earlier.
\ermrk
 \subsection{Quantum actions on Rieffel-deformed $C^*$ algebras}
We shall now consider deformation of classical manifolds and quantum group actions on them. Let $G$ be as in  Theorem \ref{main_no_action} above and assume that the rank $n$ of $G$ is at least $2$. We consider the Rieffel-deformation $C(M)_\theta$ using the left action of $T$ on $M=G/L$ (see \cite{rieffel} for the definition and details of such deformation), indexed by skew-symmetric $n \times n$ matrices $\theta$, which is a continuous field of   possibly noncommutative $C^*$ algebras. 

In a similar way, if  a compact group $K$ has an $n$-toral subgroup $T$, we can consider the Rieffel-Wang deformation  $C(K)_{\tilde{\theta}}$ (see \cite{rieffel_wang}, \cite{toral}) of 
$C(K)$ by the action of the $2n$-dimensional torus $T \times T$ on $K$ given by $(z,w)g:=zgw^{-1}$, $z,w \in T, g \in K$, and where $\tilde{\theta}=\left( \begin{array}{cc} 0 & \theta \\
-\theta & 0 \\ \end{array} \right).$ This becomes a compact quantum group with the same coalgebra structure as $C(G)$.
It is known (\cite{rieffel}) that there are canonical injective linear maps $j_\theta: C^\infty(M) \raro C(M)_\theta$ such that $C^\infty(M)_\theta:=j_\theta(C^\infty(M))$ is dense in $C(M)_\theta$ for each $\theta$.
There is a canonical action of $T$, say $\beta$, on $C(M)_\theta$, which is defined by $\beta_z(j_\theta(f))=j_\theta(L_zf),$ $L_z(f)(g):=f(zg)$, for $f \in C^\infty(M)$, $z \in T$, $g \in G$. 
\bdfn
We shall call an action $\alpha: C(M)_\theta \raro C(M)_\theta \ot \clq$ of a compact quantum group `linear' if 
 it leaves invariant the subspace $j_\theta(V) \subseteq C^\infty(M)_\theta$, where $V$ is as in Theorem \ref{main_no_action}, i.e. the subspace of $C^\infty(M)$ spanned by functions of the form $\lambda({\rm Ad}_g(H_1)),$ $\lambda \in \clg^\prime$, where $\clg$, $H_1$ are as before. 
 \edfn
\bthm
Let $\clq$ be a compact quantum quantum group with a faithful Haar state having a faithful and linear action $\alpha$ on $C(M)_\theta$. Assume furthermore that there is a quantum subgroup of $\clq$ isomorphic with the $n$-torus $T$, with the corresponding surjective morphism $\pi: \clq \raro C(T)$, such that the action $({\rm id} \ot \pi) \circ \alpha
: C(M)_\theta \raro C(M)_\theta \ot C(T)$ coincides with the canonical $T$-action $\beta$ on $C(M)_\theta$. Then $\clq$ must be isomorphic with the Rieffel-Wang deformation $C(K)_{\tilde{\theta}}$ for some compact group $K$ which acts on $M$.
\ethm
{\it Proof:}\\
Note that $(C(M)_\theta)_{-\theta}\cong C(M)$. We conclude from the proof of  Theorem 3.11 of \cite{qiso_comp} that there is an action $\alpha_{-\theta}$ of $\clq_{-\tilde{\theta}}$ on $C(M)$, which is clearly faithful and also leaves invariant $j_0(V)\equiv V$. This implies, by Theorem \ref{main_no_action}, that $\clq_{-\tilde{\theta}} \cong C(K)$ for some compact group acting on $M$. Thus, $\clq\cong C(K)_{\tilde{\theta}}$. \qed\vspace{2mm}\\

 {\bf Acknowledgement:} I would like to thank Prof. Marc A. Rieffel for inviting me to visit  the Department of Mathematics of the University of California at Berkeley, where most of the work was done. I am grateful to Prof. Rieffel for stimulating discussion with him, and also for explaining some of his work on coadjoint orbits to me, which gave me the motivation to consider the problem of quantum actions on homogeneous spaces. 
 
 I am also grateful to late Prof. S. C. Bagchi for some useful discussion on Lie groups.
 
 I would also like to thank the anonymous referee for pointing out mistakes in the older versions and constructive criticism which helped improve the paper a lot.

\end{document}